\documentclass[11pt]{article}
\usepackage{amsmath,enumerate,fancyhdr,amssymb, amsthm, pstricks, epsf, lscape}

\newcommand{\commentout}[1]{}
\newcommand{\co}[1]{}

\co{
\setlength{\oddsidemargin}{0cm} 
\setlength{\evensidemargin}{0cm}
\setlength{\textwidth}{16.0cm} 
\setlength{\topmargin}{1cm}
\setlength{\textheight}{19cm}
}

\def\eref#1{(\ref{#1})}

\newcommand{\norm}[1]{\parallel \! #1 \! \parallel}

\newcommand{\rad}[1]{\mathbb{R}^{#1}}

\newcommand{\zad}[1]{\mathbb{Z}^{#1}}
\newcommand{\qad}[1]{\mathbb{Q}^{#1}}

\newcommand{\latt}[1]{\mathbb{L}(#1)}
\newcommand{\Latt}[1]{\mathbb{L}{#1}}
\newcommand{\nlatt}[1]{\mathbb{N}(#1)}

\newcommand{\width}{{\rm width}}
\newcommand{\iwidth}{{\rm iwidth}}

\newcommand{\la}{\langle}
\newcommand{\ra}{\rangle}

\newcommand{\nin}{\noindent}

\newcommand{\tQ}{\tilde{Q}}

\newcommand{\hQ}{\hat{Q}}

\newcommand{\round}{\operatorname{round}}

\newcommand{\mydet}{\operatorname{det}}

\newcommand{\lin}{\operatorname{lin}}

\newtheorem{Claim}{Claim}

\newtheorem{Definition}{Definition}

\newtheorem{Setup}{Setup}
\newtheorem{Example}{Example}
\newtheorem{Counterexample}{Counterexample}
\newtheorem{Proposition}{Proposition}
\newtheorem{Lemma}{Lemma}
\newtheorem{Theorem}{Theorem}
\newtheorem{Corollary}[Definition]{Corollary}
\newtheorem{Remark}[Definition]{Remark}
\newtheorem{Assumption}{Assumption}
\newtheorem{Recipe}{Recipe}

\newcommand{\beq}{\begin{equation}}
\newcommand{\eeq}{\end{equation}}
\newcommand{\beqa}{\begin{eqnarray}}
\newcommand{\eeqa}{\end{eqnarray}}
\newcommand{\ba}{\begin{array}}
\newcommand{\ea}{\end{array}}
\newcommand{\bac}{\begin{array}{ccccccccccc}}
\newcommand{\eac}{\end{array}}
\newcommand{\bprop}{\begin{Proposition}}
\newcommand{\eprop}{\end{Proposition}}

\newcommand{\bcex}{\begin{Counterexample}}
\newcommand{\ecex}{\end{Counterexample}}

\newcommand{\beqast}{\begin{eqnarray*}}
\newcommand{\eeqast}{\end{eqnarray*}}
\newcommand{\benum}{\begin{enumerate}}
\newcommand{\eenum}{\end{enumerate}}
\newcommand{\bit}{\begin{itemize}}
\newcommand{\eit}{\end{itemize}}
\newcommand{\bth}{\begin{Theorem}}
\newcommand{\enth}{\end{Theorem}}

\newcommand{\bdef}{\begin{Definition}}
\newcommand{\Edef}{\end{Definition}}

\newcommand{\bsetup}{\begin{Setup}}
\newcommand{\esetup}{\end{Setup}}
\newcommand{\ble}{\begin{Lemma}}
\newcommand{\ele}{\end{Lemma}}
\newcommand{\bex}{\begin{Example}}
\newcommand{\eex}{\end{Example}}
\newcommand{\bcor}{\begin{Corollary}}
\newcommand{\ecor}{\end{Corollary}}
\newcommand{\brem}{\begin{Remark}}
\newcommand{\erem}{\end{Remark}}
\newcommand{\bass}{\begin{Assumption}}
\newcommand{\eass}{\end{Assumption}}

\newcommand{\brep}{\begin{Recipe}}
\newcommand{\erep}{\end{Recipe}}

\newcommand{\pf}[1]{\vspace{.35cm} \nin {\bf Proof {#1} }}

\newcommand{\bpx}{\begin{pmatrix}}
\newcommand{\epx}{\end{pmatrix}}

\newcommand{\bbx}{\begin{bmatrix}}
\newcommand{\ebx}{\end{bmatrix}}

\begin{document}

\title{\bf Parallel Approximation, and Integer Programming Reformulation} 
\author{G\'{a}bor Pataki and Mustafa Tural \thanks{Department of Statistics and Operations Research, UNC Chapel Hill, {\bf gabor@unc.edu, tural@email.unc.edu}} \\
Technical Report 2007-07 \\ Department of Statistics and Operations Research, UNC Chapel Hill}

\date{}

\maketitle

\begin{abstract}

\co{We show that in the knapsack feasibility problem with weight vector $a$, 
an integral vector $p$, which is short, and near parallel to $a$ gives a branching direction with small integer width. 
is a good branching direction in the sense that a small number of nodes are generated, when branching on 
$px$. Near parallelness is defined in a sense which is stronger than just requiring $|\sin(a,p)|$ to be small.
}

We show that in a knapsack feasibility problem an integral vector $p$, which is short, and near parallel to the constraint vector 
gives a branching direction with small integer width. 

We use this result to analyze two computationally efficient reformulation techniques on low density knapsack problems. 
Both reformulations have a constraint matrix with columns reduced in the sense of Lenstra, Lenstra, and Lov\'asz.
We prove an upper bound on the integer width along 
the last variable, which becomes $1, \,$ when the density is sufficiently small. 

In the proof we extract from the transformation matrices a vector which is near parallel to the constraint vector $a.$ 
The near parallel vector is a good branching direction in the original knapsack problem, and this transfers to 
the last variable in the reformulations.

\end{abstract}

\tableofcontents


\section{Introduction and notation}
\label{section-intro}

\paragraph[ip]{Geometry of Numbers and Integer Programming}

\cite{PT08}

Starting with the work of H. W. Lenstra 
\cite{L83}, algorithms based on the geometry of numbers
have been an essential part of the Integer Programming landscape. Typically, these algorithms 
reduce an IP feasibility problem to a provably small number of smaller dimensional ones, and 
have strong theoretical properties.  For instance, the algorithms of 
\cite{L83, K87, LS92} have polynomial running time in fixed dimension; the algorithm of \cite{EL05} has linear running time
in dimension two. 
One essential tool in creating the subproblems is a ``thin'' branching direction, i.e. a $c$  integral (row-)vector
with the difference between the maximum and the minimum of $cx$ over the underlying polyhedron being provably small. 
Basis reduction in lattices -- in the Lenstra, Lenstra, Lov\'asz (LLL) \cite{LLL82}, or Korkine and Zolotarev (KZ) \cite{KZ1873, K87} sense --
is usually a key ingredient in the search for a thin direction. 
For implementations, and computational results, we refer to \cite{CRSS93, GZ02, ML04}. 

A simple, and  experimentally  very successful technique for integer programming based on LLL-reduction was proposed 
by Aardal, Hurkens and A. K. Lenstra in \cite{AHL00} for equality constrained IP problems; see also
\cite{ABHLS00}.  
Consider the problem 
\beq \label{ip-eq} \tag{IP-EQ}
\ba{rcl}
Ax & = &  b \\
0 \leq & x & \leq v \\
x & \in  & \zad{n},
\ea
\eeq
where $A$ is an integral matrix with  $m$ independent rows, and let 
\beq \label{nlatt}
\nlatt{A} = \{ \, x \, \in \zad{n} \, | \,  Ax = 0 \, \}.
\eeq
The full-dimensional reformulation proposed in \cite{AHL00} is
\beq \label{ip-eq-n} \tag{IP-EQ-N}
\ba{rcl}
- x_b  \leq & V \lambda  & \leq v - x_b \\
\lambda & \in  & \zad{n-m}.
\ea
\eeq
Here $V $ and $x_b$ satisfy 
$$
\{ \, V \lambda \, | \, \lambda \in \zad{n-m} \, \} \, = \, \nlatt{A},  \, x_b \in \zad{n}, \, A x_b = b, 
$$
the columns of $V$ are reduced in the LLL-sense, and $x_b \,$ is also short. 
For several classes of hard equality constrained integer programming problems -- e.g. \cite{CD98} -- the reformulation turned out to be 
much easier to solve by commercial solvers than the original problem. 


In \cite{KP06} an experimentally just as effective reformulation method was introduced, which leaves the number of the variables the same, 
and is applicable to inequality or equality constrained problems as well.
It replaces 
\beq \label{ip} \tag{IP}
\ba{rcl}
Ax & \,\, \leq b \\
x  & \in   \zad{n} 
\ea
\eeq
with 
\beq \label{ip-r} \tag{IP-R}
\ba{rcl}
(AU)y & \,\, \leq b \\
 y    & \in   \zad{n}, & 
\ea
\eeq
where $U$ is a unimodular matrix that makes the columns of $AU$ reduced in the LLL-, or KZ-sense. 
It applies the same way, even if some of the inequalities in the IP feasibility problem are actually
equalities. Also, if the constraints are of the form $b' \leq Ax \leq b$ in \eref{ip}, the reformulation is 
just $b' \leq (AU)y \leq b, \,$  so we do not bring the system into a standard form. 
In \cite{KP06} the authors also introduced a simplified method to compute a reformulation which is essentially equivalent to \eref{ip-eq-n}.

We call \eref{ip-r} the {\em rangespace reformulation} of \eref{ip}; and \eref{ip-eq-n} the {\em nullspace reformulation} of \eref{ip-eq}.

These reformulation methods are very easy to describe (as opposed to say H. W. Lenstra's method),
but seem difficult to analyze. The only analyses are for knapsack problems, with the weight vector having a given ``decomposable'' structure, 
i.e. 
\beq \label{alambda} 
a = \lambda p + r
\eeq
with $p, r, \,$ and $\lambda \,$ integral, and $\lambda$ large with respect to $\norm{p}, \,$ and $\norm{r}$, 
see \cite{AL04, KP06}. 

The results in these papers are a first step towards a general analysis. However, besides assuming the decomposable structure a priori,
they only prove an upper bound on the width in the reformulations along the last variable. 

The goal of this paper is to prove such width results on the knapsack feasibility problem 
\beq \label{ss} \tag{KP}
\ba{rcl}
\beta_1 \, & \leq \, ax \, \leq & \, \beta_2 \\
   0 \,    & \leq \, x  \, \leq & \, v \\
           & x  \in  \zad{n}, & 
\ea
\eeq
where $a$ is a positive, integral row vector, $\beta_1, \,$ and $\beta_2 \,$ are integers without assuming any structure on $a$.
We will assume that $a$ has low density. 
The density of a set of weights $a = (a_1, \dots, a_n) \,$ is 
\beq
d(a) \, = \, \dfrac{n}{\log_2 \norm{a}_\infty}.
\eeq
Subset sum problems (when  $\beta_1 = \beta_2 = \beta, \,$ and $v$ is the vector of all ones) with the weight vector having low density 
have been extensively studied. The seminal paper of Lagarias and Odlyzko \cite{LO85} proves that 
the solution of all but at most a fraction of $1/2^n$ subset sum problems, which have a solution, and have density less than $c/n$ can be found in polynomial time, where $c \approx 4.8.$ 
Clearly $d(a) < c/n \,$ is equivalent to $2^{n^2/c} < \norm{a}_\infty$. 

Let 
\beqa
G_n(M) & = & \{ \, a \in \zad{n} \, | \, a_i \in \{ \, 1, \dots, M \, \} \}.
\eeqa
Furst and Kannan in \cite{FK89} showed that for some $c > 0$ constant, if $M \geq 2^{c n \log n}, \,$ then for almost 
all $a \in G_n(M)$ and all $\beta$ the problem \eref{ss} has a polynomial size proof of feasibility or infeasibility.
Their second result shows that for some $d > 0$ constant, if $M \geq 2^{d n^2}, \,$ then for almost 
all $a \in G_n(M)$ and all $\beta$ 
the problem \eref{ss} can be {\em solved} in polynomial time. Their proof works by constructing a candidate solution to \eref{ss}, 
and showing that for almost all $a \in G_n(M), \,$ if there is a feasible solution, then it is unique, and the candidate solution must be it. 

If we assume the availability of a {\em lattice oracle}, which finds the shortest vector in a lattice, then 
the result of \cite{LO85} can be strengthened to only requiring the density to be less than $0.6463.$ The current best result on finding the solution of 
almost all (solvable) subset sum problems using a lattice oracle is 
by Coster et al \cite{CJLOSS92}:  they require only $d(a)<0.9408.$ 

The rangespace reformulation of \eref{ss} is 
\beq \label{ss-r} \tag{KP-R}
\ba{rcl}
\beta_1 \, & \leq \, aUy \, \leq & \, \beta_2 \\
   0 \,    & \leq \, Uy  \, \leq & \, v \\
           & y \in  \zad{n}, & 
\ea
\eeq
where $U$ is a unimodular matrix that makes the columns of $\bpx a \\ I \epx U \,$ reduced in the LLL-sense (we do not analyze it with KZ-reduction).
The nullspace reformulation is 
\beq \label{ss-n} \tag{KP-N}
\ba{rcl}
- x_\beta  \leq & V \lambda  & \leq v - x_\beta \\
\lambda & \in  & \zad{n-m},
\ea
\eeq
where $x_\beta \in \zad{n}, \, a x_\beta = \beta, \, \{ \, V \lambda \, | \, \lambda \in \zad{n-m} \, \} \, = \, \nlatt{a}, \, $ and the columns of $V$ are reduced in the LLL-sense.

We will assume  $\norm{a} \, \geq \, 2^{(n/2 +1)n}. \,$ which is satisfied, when $d(a)<2/(n+2). \,$ 
We will not assume any a priori structure on $a$. In fact, a key point will be that a decomposable structure
is automatically ``discovered'' by the reformulations. Precisely, we will prove that in both reformulations a decomposition 
$a = \lambda p + r \,$ can be found from the transformation matrices, now with only $p$ integral, 
and that branching on the last variable in the reformulations will be equivalent to branching on $px$ in the original problem.

There are crucial differences between the results that {\em assume} a decomposable structure, and the results of this paper. 
For instance, in \cite{KP06} one needs to assume 
\beqa
\lambda & \geq &  2^{(n-1)/2} \norm{p} ( \norm{r}+1)^2, \\
\lambda & \geq &  2^{(n-1)/2} \norm{p}^2  \norm{r}^2,
\eeqa
for the analysis of the rangespace- and nullspace reformulations, respectively. 
A decomposition with any of these properties is unlikely to exist no matter how large $\norm{a}$ is, so we cannot 
plug the decomposition result of this paper into the argument used in \cite{KP06}.  We will  
prove a weaker lower bound on $\lambda, \,$ and an upper bound on $\norm{r}/\lambda$ in Theorems \ref{near-parallel}, and \ref{near-parallel-null},
and we will use these bounds in Theorem \ref{branch-thm} quite differently from how it is done in \cite{KP06}. 

\paragraph[not]{Notation} 
Vectors are column vectors, unless said otherwise. The $i$th unit row-vector is $e_i. \,$ 
In general, when writing $p_1, \, p_2, \,$ etc, we refer to vectors in a family of vectors. When $p_i \,$ refers to 
the $i$th component of vector $p$, we will say this explicitly. 
For a rational vector
$b \,$ we denote by 
$\round(b) \,$ the vector obtained by rounding the components of $b.$

We will assume $0 \, \leq \, \beta_1 \, \leq \, \beta_2 \, \leq \, av, \,$  and that the gcd of the components of $a$ is $1$.

For a polyhedron $Q$, and an integral row-vector $c$, the width, and the integer width of $Q$ along $c$ are
\beqast
\width(c, Q) & = & \max \, \{ \, cx \, | \, x \in Q \, \}  - \min \, \{ \, cx \, | \, x \in Q \, \}, \; \text{and} \\
\iwidth(c, Q) & = & \lfloor \max \, \{ \, cx \, | \, x \in Q \, \} \rfloor - \lceil \min \, \{ \, cx \, | \, x \in Q \, \} \rceil + 1.
\eeqast
The integer width is 
the number of nodes generated by branch-and-bound when branching on the hyperplane $cx$; in particular,
$\iwidth(e_i,Q)$ is the number of nodes generated when branching on $x_i$.
If the integer width along any integral vector is zero, then $Q$ has no integral points.
Given an integer program labeled by ${\rm (P)}, \,$ and $c \,$ an integral vector, we also write 
$\width(c, {\rm (P)}),$ and $\iwidth(c, {\rm (P)})$ for the width, and the integer width of the LP-relaxation of ${\rm (P)}$ along $c, $ respectively.

A lattice in $\rad{n}$ is a set of the form 
\beq \label{def-latt-B}
L \, = \, \latt{B} \, = \, \{ \, Bx \, | \, x \in \zad{n} \, \},
\eeq
where $B$ is a real matrix with $n$ independent columns, called a {\em basis} of $L$. 
A square, integral matrix $U$ is {\em unimodular} if $\det U = \pm 1$. It is well known that 
if $B_1$ and $B_2$ are bases of the same lattice, then $B_2 = B_1 U \,$ for some unimodular $U$. The determinant of $L \,$ is 
\beq
\mydet L \, = \, (\det B^T B)^{1/2},
\eeq
where $B$ is a basis of $L$; it is easy to see that $\mydet L$ is well-defined.

The LLL basis reduction algorithm \cite{LLL82} computes a reduced basis of a lattice in which
the columns are ``short'' and ``nearly'' orthogonal. 
It runs in polynomial time for rational lattices. For simplicity, we use Schrijver's definition from  \cite{S86}.
Suppose that $B$ has $n$ independent columns, i.e. 
\beq
B = [ b_1, \dots, b_n ], \; 
\eeq
and $b_1^*, \dots, b_n^*$ form the Gram-Schmidt orthogonalization of $b_1, \dots, b_n, \,$ that is $b_1 = b_1^*, \,$ and 
\beq \label{bibist}
b_i =  b_i^* + \sum_{j=1}^{i-1} \mu_{ij} b_j^* \,\, \text{with } \,\, \mu_{ij} \, = \, b_i^T b_j^* /\norm{b_j^*}^2 \,\, (i=2, \dots, n; \, j \leq i-1).
\eeq
We call $b_1, \dots, b_n$ an {\em LLL-reduced basis of } $\latt{B}, \,$ if 
 \beqa \label{mucond}
| \mu_{ij} | & \leq & 1/2 \; \,\,\, (i=2, \dots, n; \, j = 1, \dots, i-1), \, \text{and} \\ \label{exch-cond}
\norm{b_i^*}^2 & \leq &  2 \norm{b_{i+1}^*}^2 \,\, (i=1, \dots, n-1).
\eeqa
For an integral lattice $L, \,$ its {\em orthogonal lattice} is defined 
as
$$
L^\perp \, = \, \{ \, y \in \zad{n} \, | \, y^T x  = 0 \; \forall x \in L \, \},
$$
and it holds that (see e.g. \cite{M03}) 
\beq \label{orthlattice}
\mydet L^\perp \, \leq \, \mydet L.
\eeq
Suppose $A$ is an integral matrix with independent rows. Then recalling 
\eref{nlatt}, $\nlatt{A} \, $ is the same as 
$\latt{A^T}^\perp. \,$ 
A lattice $L \subseteq \zad{n}$ is called {\em complete}, if 
$$
L \, = \, \lin \, L \, \cap \, \zad{n}.
$$
The following lemma summarizes some basic results in lattice theory that we will use later on;
for a proof, see for instance \cite{M03}.
\ble \label{compl-lattice-lemma}
Let $V$ be an integral matrix with $n$ rows, and $k$  independent columns, and $L = \latt{V}$. 
Then (1) through (3) below are equivalent.
\benum
\item \label{compl-lattice-lemma-1} $L$ is complete;
\item \label{compl-lattice-lemma-2} $\mydet L^\perp  = \mydet L$;
\item \label{compl-lattice-lemma-3} There is a unimodular matrix $Z$ s.t. 
$$
Z V \, = \, \bpx I_k \\ 0_{(n-k) \times k} \epx.
$$
\eenum
Furthermore, if $Z$ is as in part \eref{compl-lattice-lemma-3}, then the last $n-k$ rows of $Z$ are a basis of  $L^\perp$.
\ele
\qed

For an $n$-vector $a, \,$ we will write 
\beq
\ba{rcl}
f(a) & = & 2^{n/4}/\norm{a}^{1/n} \\
g(a) & = & 2^{(n-2)/4}/\norm{a}^{1/(n-1)}.
\ea
\eeq

\section{Main results}
\label{section-main-results}

In this section we will review the main results of the paper, give some examples, explanations, and some proofs that show their connection. The bulk of the work 
is the proof of Theorems \ref{near-parallel}, \ref{near-parallel-null}, and \ref{branch-thm}, which is done in Section \ref{proofs-section}. 

The main purpose of this paper is an analysis of the reformulation methods. This is done in Theorem 
\ref{main}, which proves an upper bound on the number of branch-and-bound nodes, when branching on the last variable in the reformulations.
However, some of the intermediate results may be of interest on their own right. 

We will rely on Theorem \ref{sublat-thm}, proven in the companion paper \cite{PT08}, 
which gives a bound on the determinant of a sublattice in an LLL-reduced basis, thus generalizing 
the well-known result from \cite{LLL82} showing that the first vector in such a basis is short. 

Theorems \ref{near-parallel} and \ref{near-parallel-null} show that an integral vector $p, \,$ which is ``near parallel'' to $a$ can be extracted from the transformation 
matrices of the reformulations. The notion of near parallelness that we use is stronger than 
just requiring $|\sin(a,p)| \,$ to be small, and the relationship of the two parallelness concepts is clarified in Proposition \ref{rl-prop}. 
A method to find a near parallel vector using simultaneous diophantine approximation was described by Frank and Tardos 
in \cite{FT87}. Their goal was quite different from ours, and a near parallel vector derived via diophantine approximation 
is not suitable for the analysis of the reformulation methods.
For completeness, we will give an overview of their method in subsection \ref{sub-da}. 

Theorem \ref{branch-thm} proves an upper bound on $\iwidth(p, \eref{ss}), \,$ where $p$ is an integral vector. 
A novelty of the bound is that it does not depend on $\beta_1, \,$ and $\beta_2, \,$ only on their difference. 
We show through examples that 
this bound is quite useful when $p \,$ is a near parallel vector found according to Theorems \ref{near-parallel} and \ref{near-parallel-null}.

In the end, a transference result between branching directions in the original, and reformulated problems completes the proof of 
Theorem \ref{main}.

\bth \label{main} Suppose $\norm{a} \, \geq \, 2^{(n/2 +1)n}. \,$ Then 
\benum
\item \label{main1}
$
\iwidth{(e_n, \text{ {\rm \eref{ss-r}} } \!\!)}  \, \leq \, \lfloor \, f(a) ( 2 \norm{v} + (\beta_2 - \beta_1) ) \rfloor + 1.
$  
\item \label{main2}
$ 
\iwidth{(e_{n-1}, \text{ {\rm \eref{ss-n}} } \!\!)}  \, \leq \,\lfloor 2 g(a) \norm{v} \rfloor + 1.
$  
\eenum
\enth
\qed

The integer width, and the width differ by at most one, 
and are frequently used interchangeably in integer programming 
algorithms. For instance, the algorithms of \cite{L83, LS92}  find a branching direction in which the width is bounded by an exponential function 
of the dimension. The goal is proving polynomial running time in fixed dimension, and this would still be achieved 
if the width were larger by a constant. 

In contrast, when $\norm{a} \,$ is sufficiently large, Theorem \ref{main} implies that 
the integer width is at most {\em one} in both reformulations. 


The following was proven in \cite{PT08}: 
\bth
\label{sublat-thm}
Suppose that $b_1, \dots, b_n \,$ form an LLL-reduced basis of the lattice $L, \,$ and denote by 
$L_\ell \,$ the lattice generated by $b_1, \dots, b_\ell. \,$ 
Then 
\beq \label{dldn} 
\mydet L_\ell \, \leq \, 2^{\ell(n-\ell)/4} (\mydet L)^{\ell/n}.
\eeq
\enth
\nin Theorem \ref{sublat-thm} is a natural generalization of 
$\norm{b_1} \leq 2^{(n-1)/4} (\mydet L)^{1/n}$ (see \cite{LLL82}). 

Given $a$ and $p$ integral vectors, we will need the notion of their near parallelness. The obvious thing would be 
to require that $|\sin(a,p)|$ is small. Instead, we will write a decomposition 
\beq \label{decomp0} \tag{DECOMP}  
a = \lambda p + r, \, \text{with} \, \lambda \in \qad{}, \, r \in \qad{n}, \, r \bot p,
\eeq
and ask for $\norm{r}/\lambda$ to be small. The following proposition clarifies the connection of the two near parallelness concepts, 
and shows two useful consequences of the latter one.
\bprop \label{rl-prop}
Suppose that $a, p \, \in  \, \zad{n}, \,$ and $r$ and $\lambda$ are defined to satisfy \eref{decomp}. Assume w.l.o.g. $\lambda > 0.$ 
Then 
\benum
\item \label{rl-prop1} $\sin(a,p)  \, \leq \, \norm{r}\!/\lambda.$ 
\item \label{rl-prop2} For any $M \,$ there is $a, \,p \,$ with $\norm{a} \geq M \,$ such that the inequality in \eref{rl-prop1} is strict.
\item \label{rl-prop3} Denote by $p_i \,$ and $a_i \,$ the $i$th component of $p, \,$ and $a. $
If $\norm{r}/\lambda < 1, \,$ and $p_i \neq 0, \,$ then the signs of $p_i \,$ and $a_i \,$ agree. Also, 
if $\norm{r}\!/\lambda < 1/2,$ then  $\lfloor a_i/\lambda \rceil = p_i.$ 
\eenum
\eprop
\pf{} Statement  \eref{rl-prop1} follows from 
\beq \label{sinapleq}
\ba{rclclcl}
\sin(a,p) & = &  \norm{r}/\norm{a} \, \leq \,  \norm{r}/\norm{ \lambda p} \, \leq \, \norm{r}/\lambda,
\ea
\eeq
where in the last inequality we used the integrality of  $p.$ 

To see \eref{rl-prop2}, one can choose $a$ and $p$ to be near orthogonal, to make $\norm{r}/\lambda$ arbitrarily large, while $\sin(a,p)$ will always be bounded 
by $1$. A more interesting example is from considering the family of $a, \,$ and $p \,$ vectors
\beq
\ba{rcl}
a & = & \bpx m^2+1, & m^2 \epx, \\
p & = & \bpx m+1, & m \epx
\ea
\eeq
with $m$ an integer. Letting $\lambda \,$ and $r \,$ be defined as in the statement of the proposition, a straightforward computation (or experimentation) shows 
that as $m \rightarrow \infty$ 
\beqast
\sin(a, p) & \rightarrow & 0, \,\,\, \\ 
\norm{r}/\lambda & \rightarrow & 1/\sqrt{2}. 
\eeqast
Statement \eref{rl-prop3} is straighforward from 
\beq
\ba{rcl} \label{earlier}
a_i/\lambda & = & p_i +  r_i/\lambda. 
\end{array}
\eeq
\qed

The next two theorems show how the near parallel vectors can be found from the transformation matrices of the reformulations.
\bth \label{near-parallel}
\label{pa-1}
Suppose $\norm{a} \, \geq \, 2^{(n/2 +1)n}. \,$ Let $U$ be a unimodular matrix such that 
the columns of 
$$
\bpx
a \\
I
\epx U 
$$
are LLL-reduced, and $p$ the last row of $U^{-1}$. 
Define $r$ and $\lambda$ to satisfy \eref{decomp}, and assume w.l.o.g. $\lambda > 0.$ 

\nin Then 
\benum
\item \label{pa-1-1} $\norm{p} (1+\norm{r}^{2})^{1/2} \leq \norm{a} f(a)$;
\item \label{pa-1-2} $\lambda \geq 1/f(a)$; 
\item \label{pa-1-3} $\norm{r}/\lambda \leq 2 f(a)$.
\eenum
\enth
\qed

\bth \label{near-parallel-null}
Suppose $\norm{a} \, \geq \, 2^{(n/2 +1)n}. \,$ Let 
$V$ be a matrix whose columns are an LLL-reduced basis of $\nlatt{a}$, $b$ an integral column vector with $ab=1$, 
and $p$ the $(n-1)st$ row of $(V,b)^{-1}$. 
Define $r$ and $\lambda$ to satisfy \eref{decomp}, and assume w.l.o.g. $\lambda > 0.$ 

Then $r \neq 0, \,$ and 
\benum
\item \label{pa-null-1-1} $\norm{p} \norm{r} \leq \norm{a} g(a)$;
\item \label{pa-null-1-2} $\norm{r}/\lambda \leq 2 g(a)$.
\eenum
\enth
\qed

It is important to note that $p \,$ is integral, but $\lambda$ and $r \,$ may  not be.
Also, the measure of parallelness to $a, \, $ i.e. the upper bound on $\norm{r}/\lambda \,$ is quite similar 
for the $p$ vectors found in Theorems  \ref{near-parallel} and \ref{near-parallel-null}, but their length can be 
quite different. When $\norm{a}$ is large, the $p $ vector in Theorem \ref{near-parallel} is guaranteed to be much shorter than $a \,$ by 
$\lambda \geq 1/f(a). \,$ On the other hand, the $p$ vector from Theorem \ref{near-parallel-null} may be much {\em longer} than $a: $ 
the upper bound on $\norm{p} \norm{r} \,$ does not guarantee any bound on $\norm{p}, \,$ since $r \,$ can be fractional.

The following example illustrates this: 
\bex \label{ex1}
Consider the vector 
\beq
\ba{rcl}
a & = & \bpx 3488, &    451, &    1231, &    6415, &   2191 \epx.
\ea
\eeq
We computed $p_1, \, r_1, \, \lambda_1 \,$  according to Theorem \ref{near-parallel}:
\beq
\ba{rcl}
p_1 & = & \bpx  62, &   8, &  22, &  114, &  39 \epx, \\
r_1 & = & \bpx 0.2582, &    0.9688, &  -6.5858, &  2.0554, &   -2.9021 \epx, \\
\lambda_1 & = & 56.2539, \\
\norm{r_1}/\lambda_1 & = & 0.1342.
\ea
\eeq
We also computed $p_2, \, r_2, \, \lambda_2 \,$  according to Theorem \ref{near-parallel-null}; note 
$\norm{p_2} > \norm{a}$:
\beq
\ba{rcl}
p_2 & = & \bpx 12204, &   1578, &   4307, &  22445, &   7666 \epx \\
r_2 & = & \bpx -0.0165, & -0.0071, &  0.0194, &   0.0105, &   -0.0140 \epx \\
\lambda_2 & = &  0.2858 \\
\norm{r_2}/\lambda_2 & = &  0.1110.
\ea
\eeq
\eex
\qed

Theorem \ref{branch-thm} below gives an upper bound on the number of branch-and-bound nodes when branching on a hyperplane 
in \eref{ss}.
\bth \label{branch-thm}
Suppose that $a = \lambda p + r, \,$ with $p \geq 0.$ Then 
\beq
\iwidth(p, \eref{ss}) \, \leq \, \left\lfloor \dfrac{\norm{r} \norm{v}}{\lambda} + \dfrac{\beta_2 - \beta_1}{\lambda} \right\rfloor + 1.
\eeq
\co{
\beq
\iwidth(p, \eref{ss}) \, \leq \, (\beta_2 - \beta_1)/\lambda + \norm{r} \norm{v}/\lambda + 1.
\eeq
}
\enth
This bound is quite strong for 
near parallel vectors computed from Theorems \ref{near-parallel} and \ref{near-parallel-null}. For instance, let  $a, \, p_1, \, r_1, \lambda_1 \,$ 
be as in Example \ref{ex1}. If $\beta_1 = \beta_2 \,$ in a knapsack problem with weight vector $a, $ 
and each $x_i$ is bounded between $0$ and $11, \,$ then Theorem \ref{branch-thm} implies that the 
integer width is at most one. 
At the other extreme, it also implies that the integer width is at most one,
if each $x_i$ is bounded between $0$ and $1$, and $\beta_2 - \beta_1 \leq 39. \,$ 
However, this bound does not seem as useful, when $p$ is a ``simple'' vector, say a unit vector.

We now complete the proof of Theorem \ref{main}, based on a simple transference result between branching directions,
taken from \cite{KP06}. 

\pf{of Theorem \ref{main}} 

Let us denote by $Q, \, \tQ, \,$ and $\hQ \,$ the feasible sets of the LP-relaxations of 
\eref{ss}, of \eref{ss-r}, and of \eref{ss-n}, respectively. 

First, let $U, \,$ and $p$ be the transformation matrix, and the near parallel vector from  Theorem \ref{near-parallel}. 
It was shown in \cite{KP06} that $\iwidth{ (p, Q) } \, = \, \iwidth{ (pU, \tQ) }.$
\co{\beq \label{enp}
\ba{rcl}
\iwidth{ (p, Q) } & = & \iwidth{ (pU, \tQ) }. 
\ea
\eeq}
But  $pU = \pm e_n, \,$ so 
\beq \label{enp}
\ba{rcl}
\iwidth{ (p, Q) } & = & \iwidth{ (e_n, \tQ) }. 
\ea
\eeq
On the other hand, 
\beq \label{iw-range}
\ba{rcl} 
\iwidth{ (p, Q) } & \leq & \left\lfloor \dfrac{\norm{r} \norm{v}}{\lambda} + \dfrac{\beta_2 - \beta_1}{\lambda} \right\rfloor + 1 \\
                  & \leq & \lfloor \, f(a) ( 2 \norm{v} + (\beta_2 - \beta_1) ) \rfloor + 1
\ea
\eeq
with the first inequality coming from Theorem \ref{branch-thm}, and the second from 
using the bounds on $1/\lambda \,$ and $\norm{r}/\lambda \,$ from Theorem \ref{near-parallel}.
Combining \eref{enp} and \eref{iw-range} yields \eref{main1} in Theorem \ref{main}.

Now let $V, \,$ and $p$ be the transformation matrix, and the near parallel vector from  Theorem \ref{near-parallel-null}. 
It was shown in \cite{KP06} that $\iwidth{ (p, Q) } \, = \, \iwidth{ (pV, \hQ) }.$ But $pV = \pm e_{n-1}, \,$ so 
\beq \label{enp-null}
\ba{rcl}
\iwidth{ (e_{n-1}, \hQ) } & = & \iwidth{ (p, Q) }.
\ea
\eeq
On the other hand, 
\beq \label{iw-null}
\ba{rcl} 
\iwidth{ (p, Q) } & \leq & \left\lfloor \dfrac{\norm{r} \norm{v}}{\lambda} \right\rfloor + 1 \\
                  & \leq & \lfloor \, g(a) ( 2 \norm{v} ) \rfloor + 1.
\ea
\eeq
with the first inequality coming from Theorem \ref{branch-thm}, and the second from 
using the bound on $\norm{r}/\lambda \,$ in Theorem \ref{near-parallel-null}.
Combining \eref{enp-null} and \eref{iw-null} yields \eref{main2} in Theorem \ref{main}.

\co{Since 
$$
\dfrac{39}{\lambda_1} + \dfrac{\iwidth{ (p, Q) } & = & \iwidth{ (pU, \tQ) }. \norm{r_1}{\lambda_1}} \sqrt{5} < 1,
$$
as long as the }

\section{Proofs}
\label{proofs-section}

\subsection{Near parallel vectors: intuition, and proofs for Theorems \ref{near-parallel} and \ref{near-parallel-null}}

\nin{\bf Intuition for Theorem \ref{near-parallel}} We review a proof from \cite{KP06}, which applies 
when we know {\em a priori} the existence of a decomposition
\beq \label{decomp}
a = p \lambda + r,
\eeq
with $\lambda$ large with respect to $\norm{p}, \,$ and $\norm{r}. \,$ The reason that the columns of 
$$
\bpx
a \\
I
\epx \, = \, \bpx
\lambda p + r \\
I
\epx
$$
are {\em not} short and orthogonal is the presence of the $\lambda_i p_i \,$ components in the first row. 
So if postmultiplying by a unimodular $U$ results in reducedness, it is natural to expect that many components of
$pU \,$ will be zero; indeed it follows from the properties of LLL-reduction, that the first $n-1$ components {\em will} be zero. Since $U \,$ has full rank, 
the $n$th component of $pU \,$ must be nonzero. So $p \,$ will be the a multiple of the last row of $U^{-1}, \,$ in other words, the last row of 
$U^{-1} \,$ will be near parallel to $a. \,$ 
(In \cite{KP06} it was assumed that $p, \, r, \,$ and $\lambda \,$ are integral, but the proof would work even if $\lambda \,$ and $r \,$ were rational. )

It is then natural to expect that the last row of $U^{-1} \,$ will give a near parallel vector to $a, \,$ even if a decomposition like 
\eref{decomp} is not known in advance. This is indeed what we show in Theorem \ref{near-parallel}, when $\norm{a} \,$ is sufficiently large.

\pf{ of Theorem \ref{near-parallel}} First note that the lower bound on $\norm{a}$ implies
\beq \label{34}
f(a) \leq \sqrt{3}/2.
\eeq
Let 
$L_\ell$ be the lattice generated by the first $\ell$ columns 
of $\bpx a  \\ I \! \epx U, $ and 
\beqast
Z & = & \bpx 0 & U^{-1} \\
             1 & - a  
        \epx.
\eeqast
Clearly, $Z$ is unimodular, and 
\beqa \label{zau}
Z \bpx aU  \\ U \epx & = & \bpx I_{n} \\  0_{1 \times n} \epx.
\eeqa
So Lemma \ref{compl-lattice-lemma} implies that 
$L_\ell$ is complete, and the last $n+1-\ell$ rows of $Z$ generate 
$L_\ell^\perp$. The last row of $Z$ is $(1, -a), \,$ and the next-to-last is $(0, p), \,$ so we get
\beq \label{lnm1-1}
\ba{rcl} 
\mydet L_{n} & = & \mydet L_{n}^\perp \, = \, (\norm{a}^2 + 1)^{1/2},  \\ 
\mydet L_{n-1} & = & \mydet L_{n-1}^\perp \, = \, \norm{p} (1+\norm{r}^{2})^{1/2}. 
\end{array}
\eeq
Theorem \ref{sublat-thm} implies 
\beq \label{mydet}
\ba{rcl}
\mydet \, L_{n-1} & \leq & 2^{(n-1)/4} (\mydet L_n)^{1 - 1/n}.
\end{array}
\eeq
Substituting into \eref{mydet} from \eref{lnm1-1} gives 
\beq
\ba{rcl}
\norm{p}(1+\norm{r}^2)^{1/2} & \leq & 2^{(n-1)/4} (\sqrt{\norm{a}^2 + 1})^{1 - 1/n}  \\
                             & \leq & 2^{n/4} \norm{a}^{1-1/n} \\
                             & = & \norm{a} f(a),
\end{array}
\eeq
with the second  inequality coming the lower bound on $\norm{a}$. This shows \eref{pa-1-1}.

\pf{ of \eref{pa-1-2}} From \eref{pa-1-1} we directly obtain
\co{\beq \label{rlambda}
\dfrac{f(a)^2 \norm{a}^2 - \norm{r}^2}{\norm{p}^2} \, \geq \, 1 \, = \, \dfrac{f(a)^2 \norm{a}^2}{f(a)^2 \norm{a}^2}.
\eeq}
\beq \label{rlambda}
\ba{rcl}
\dfrac{f(a)^2 \norm{a}^2 - \norm{r}^2}{\norm{p}^2} & \geq & \dfrac{f(a)^2 \norm{a}^2 - \norm{p}^2 \norm{r}^2}{\norm{p}^2} \\
                                                              & \geq &  1 \\
                                                              & = & \dfrac{f(a)^2 \norm{a}^2}{f(a)^2 \norm{a}^2},
\end{array}
\eeq
where in the first inequality we used $\norm{p} \geq 1$. Now note 
$$
\norm{p}^2 \leq f(a)^2 \norm{a}^2, \,
$$ 
i.e. the the denominator of the first expression in \eref{rlambda} is not larger than the denominator of the
last expression. So if we replace $f(a)^2$ by $1$ in the {\em numerator} of both, the inequality will remain valid. The result is
\beq \label{rlambda-1}
\dfrac{\norm{a}^2 - \norm{r}^2}{\norm{p}^2} \, \geq \, \, \dfrac{1}{f(a)^2},
\eeq
which is the square of the required inequality.

\pf{ of \eref{pa-1-3}} We have 
\beq
\ba{rcl}
\dfrac{\norm{r}^2 }{\lambda^2} & \leq &  \dfrac{\norm{p}^2 \norm{r}^2 }{ \norm{\lambda p}^2} \\
                               & \leq & \dfrac{\norm{p}^2 \norm{r}^2 }{\norm{a}^2 - \norm{r}^2} \\
                               & \leq & \dfrac{f(a)^2 \norm{a}^2}{\norm{a}^2 - \norm{r}^2} \\
                               & \leq & \dfrac{f(a)^2 \norm{a}^2}{\norm{a}^2 - f(a)^2 \norm{a}^2} \\
                               & =    & \dfrac{f(a)^2 }{1  - f(a)^2} \\
                               & \leq   & 4 f(a)^2,
\ea
\eeq
where the first inequality comes from Proposition \ref{rl-prop}, the last from \eref{34}, and the others 
are straightforward. 

\qed

\nin{\bf Intuition for Theorem \ref{near-parallel-null}} We recall  a proof from \cite{KP06}, which applies 
when we know {\em a priori} the existence of a decomposition like in \eref{decomp} with 
$\lambda$ large with respect to $\norm{p}, \,$ and $\norm{r}, \,$ and $p $ not a multiple of $r. $ 
It is shown there that the first $n-2$ components of $pV$ will be zero. 
Denote by $L_{\ell}$ the lattice generated by the first $\ell$ columns of $V$. 
So $p$ is in $L_{n-2}^\perp, \,$ and it is not a multiple of $a, \,$ but it is near parallel to it.

So one can  expect that an element of $L_{n-2}^\perp \,$ which is distinct from $a \,$ will be near parallel to $a, $
even if a decomposition like \eref{decomp} is not known in advance. The $p$ described in Theorem 
\ref{near-parallel-null}  will be  such a vector.


\pf{of Theorem \ref{near-parallel-null}} The lower bound on $\norm{a}$ implies
\beq \label{34-null}
g(a) \leq \sqrt{3}/2.
\eeq
As noted above, let 
$L_\ell$ be the lattice generated by the first $\ell$ columns 
of $V.$ 
We have 
\beqa \label{vb}
(V,b)^{-1} V & = & \bpx I_{n-1}  \\ 0  \epx.
\eeqa
So Lemma \ref{compl-lattice-lemma} implies that 
$L_\ell$ is complete, and the last $n-\ell$ rows of $(V,b)^{-1}$ generate 
$L_\ell^\perp$. It is elementary to see that the last row of $(V,b)^{-1}$ is $a,$ and by definition the next-to-last row is $p,$ and 
these rows are independent, so $r \neq 0. $ Also, 
\beq \label{lnm1-2}
\ba{rcl} 
\mydet L_{n-1} & = & \mydet L_{n-1}^\perp \, = \, \norm{a},  \\ 
\mydet L_{n-2} & = & \mydet L_{n-2}^\perp \, = \, \norm{p} \norm{r}. 
\end{array}
\eeq
Theorem \ref{sublat-thm} with $n-1$ in place of $n, \,$ and $n-2$ in place of $\ell$ implies 
\beq \label{mydet-2}
\ba{rcl}
\mydet \, L_{n-2} & \leq & 2^{(n-2)/4} (\mydet L_{n-1})^{1 - 1/(n-1)}.
\end{array}
\eeq
Substituting into \eref{mydet-2} from \eref{lnm1-2} gives 
\beq
\ba{rcl}
\norm{p}\norm{r} & \leq & 2^{(n-2)/4} \norm{a}^{1-1/(n-1)} \\
                             & = & \norm{a} g(a),
\end{array}
\eeq
as required. 

\pf{of \eref{pa-null-1-2}} It is enough to note that in proof of \eref{pa-1-3} in Theorem \ref{near-parallel} we only used the 
inequality $\norm{p}^2 \norm{r}^2 \leq f(a)^2 \norm{a}^2. \,$ So the exact same argument works here as well with $g(a)$ instead of $f(a), $
and invoking \eref{34-null} as well. 

\qed

\subsection{Branching on a near parallel vector: proof of Theorem \ref{branch-thm}}

This proof is somewhat  technical, so we state, and prove some intermediate claims, to improve readability. 
Let us fix $a, \, p, \, \beta_1, \, \beta_2, $ and $v. \,$
For a row-vector $w, \,$ and an integer $\ell \,$ we write
\beq
\ba{rcl}
\max(w,\ell) & = & \max \, \{ \, wx \, | \, px \leq \ell, \, 0 \leq x \leq v \, \} \\
\min(w,\ell) & = & \min \, \{ \, wx \, | \, px \geq \ell, \, 0 \leq x \leq v \, \}.
\end{array}
\eeq
The dependence on $p, \,$ on $v, \,$ and on the sense of the constraint (i.e. 
$\leq, \,$ or $\geq \,$) is not shown  by this notation; however,
we always use $px \leq \ell \,$ with ``max'', and $px \geq \ell \,$ with ``min'', and $p \, $ and $v$ are fixed.
Note that as $a$ is a row-vector, and $v$ a column-vector, $av$ is their inner product, and the meaning of $pv$ is similar.

\begin{Claim} \label{pe-claim} Suppose that $\ell_1$ and $ \ell_2$ are integers in $\{ 0, \dots,  pv \}.$ Then 
\beqa \label{minmax}
\min(a, \ell_2) - \max(a,\ell_1) & \geq & - \norm{r} \norm{v} + \lambda(\ell_2-\ell_1). 
\eeqa
\end{Claim}
\pf{} The decomposition of $a$ shows
\beq \label{max-and-min}
\ba{rcl}
\max(a,\ell_1)    & \leq & \max(r,\ell_1) + \lambda \ell_1,  \; \text{and} \\ 
\min(a,\ell_2)    & \geq & \min(r,\ell_2) + \lambda \ell_2.
\end{array}
\eeq
So we get the following chain of inequalities, with ensuing explanation: 
\beq \label{minmax-l12}
\ba{rcl}
\min(a, \ell_2) - \max(a,\ell_1) & \geq & \min(r, \ell_2) - \max(r,\ell_1) + \lambda(\ell_2 - \ell_1) \\
                         & \geq &             r x_2 - r x_1 + \lambda(\ell_2 - \ell_1) \\
                         & = &   r (x_2 - x_1) + \lambda(\ell_2 - \ell_1) \\
                         & \geq & - \norm{r} \norm{v} + \lambda(\ell_2 - \ell_1).
\end{array}
\eeq
Here $x_2$ and $x_1$ are the solutions that attain the maximum, and the minimum in $\min(r, \ell_2)$ and 
$\max(r,\ell_1), \,$ respectively. The last inequality follows from the fact that the $i$th component of $x_2 - x_1 \, $ is at most $v_i$ in absolute value,
and the Cauchy-Schwartz inequality.

\nin{\bf End of proof of Claim \ref{pe-claim}}

\nin Next, let us note
\beqa \label{minmax1}
\min(a,k)  & \leq  & \max(a,k) \,\, \text{for} \,\, k \in \{0, \dots, pv \}.
\eeqa
Indeed, \eref{minmax1} holds, since 
the feasible sets of the optimization problems defining $\min(a,k), \,$ and $\max(a,k) \,$ contain $\{ \, x \, | \, px=k, \, 0 \leq x \leq v \, \}. \,$ 

\nin The nonnegativity of $p \,$ and of $a \,$ imply
$
\min(a,0) = 0, \, \text{and} \, \max(a, pe) = av.
$
The proof of the following claim is trivial, hence omitted.
\begin{Claim} \label{branch-claim} Suppose that $\ell_1$ and $ \ell_2$ are integers in $\{ 0, \dots,  pv \}$ with $\ell_1 + 1 \leq \ell_2, \,$ and 
\beq \label{maxminl1l2}
\max(a, \ell_1) < \beta_1 \leq \beta_2  < \min(a, \ell_2).
\eeq
Then for all $x \,$ with $\beta_1 \leq ax \leq \beta_2, \, 0 \leq x \leq v$ 
\beq \label{minl1}
\ell_1 < px < \ell_2
\eeq
holds.
\end{Claim}

\nin We assume for simplicity 
\beq \label{maxa0beta}
\max(a,0) < \beta_1 \leq \beta_2 < \min(a,pe);
\eeq
the cases when this fails to hold are easy to handle separately. 
Let $\ell_1$  be the largest, and $\ell_2$ the smallest integer such that 
\beq
\max(a,\ell_1) < \beta_1 \leq \beta_2 < \min(a,\ell_2).
\eeq
From \eref{minmax1} $\ell_2 \geq \ell_1 + 1 \,$ follows, and 
Claim \ref{branch-claim} yields
\beq \label{iwidthclaim}
\iwidth{ (p, \text{ {\rm \eref{ss}} } \!\!) }  \, \leq \, \ell_2 - \ell_1 - 1.
\eeq
By the choices of $\ell_1,  \,$ and $\ell_2 \,$ we have 
\beq
\beta_1 \leq \max(a, \ell_1 +1), \,\, \text{and} \,\, \beta_2 \geq \min(a, \ell_2 -1),
\eeq
hence Claim \ref{pe-claim} leads to 
\beq \label{b1b2}
\ba{rcl}
\beta_2 - \beta_1 & \geq & \min(a, \ell_2 -1) - \max(a, \ell_1 +1) \\
                  & \geq & - \norm{r} \norm{v} + \lambda(\ell_2 - \ell_1 -2),
\ea
\eeq
that is
\beq \label{b1b2-2}
\ba{rcl}
\ell_2 - \ell_1 -2 & \leq & \dfrac{\beta_2 - \beta_1}{\lambda} + \dfrac{\norm{r} \norm{v}}{\lambda}.
\end{array}
\eeq
Comparing \eref{iwidthclaim} and \eref{b1b2-2} yields completes the proof.

\qed

\section{Discussion} 

\subsection{Connection with diophantine approximation, and other notions of near parallelness}
\label{sub-da}

Given a rational vector $b, \,$ simultaneous diophantine approximation (see e.g. \cite{LLL82, L85}) computes an integral  vector $p, \,$ and an integer 
$q, \,$ such that $q, \,$ and $\norm{b - (1/q)p} \,$ are both small. Frank and Tardos in \cite{FT87} has explored the following methodology to compute 
a vector $p$ that is near parallel to an {\em integral} vector $a. \,$ 
They apply diophantine approximation to $(1/\norm{a}_\infty a, \,$  then set 
$\lambda = \norm{a}_\infty/q, \, r = a - \lambda p. \,$ Then $\norm{r}/\lambda \,$ will be small, and if $\norm{a} \,$ is large, then $\lambda$ will be large.
\footnote{Thanks are due to Laci Lov\'asz and Fritz Eisenbrand for pointing out this connection}.

The relevance of Theorems \ref{near-parallel} and \ref{near-parallel-null} is not just finding near parallel vectors: 
it is finding  a near parallel $p, $ which corresponds to a unit vector in the rangespace- and nullspace reformulations, thus leading to the analysis 
of Theorem \ref{main}. 

Finding an integral vector, which is  near parallel to an other integral or rational one  has other applications as well.
In \cite{HN04} Huyer, and Neumaier studied several notions of near parallelness,  
presented numerical algorithms, and applications to verifying the feasibility of a linear system of inequalities.

\subsection{Successive approximation} 

Theorems \ref{near-parallel} and \ref{near-parallel-null} approximate $a$ by a single vector.
It is natural to ask: if one row of $U^{-1}$, or of 
$(V, b)^{-1}$ is a good approximation of $a$, can we construct a better 
approximation from $2, 3, \dots, k \,$ rows? 

The answer is yes, and we outline the corresponding results below, and their proofs, which are 
slight modifications of the proofs of Theorems \ref{near-parallel} and \ref{near-parallel-null}. 
As of now, we don't know how to use the general results for a 
better analysis of the reformulations than what is already given in Theorem \ref{main}.  

So we mainly state the successive approximation results for the interesting geometric intuition they give.
\co{
However, there is a natural 
geometric intuition behind them: if one row of $U^{-1} \,$ or of $(V,b)^{-1} \,$ gives a good
approximation of $a, \,$  then a combination of $2, 3, \dots, k \,$ must give increasingly better approximations. 
Since Theorems \ref{near-parallel-k} and \ref{near-parallel-null-k} verify this intuition, it 
is worth stating them, and outlining the proofs.
}
\co{
We state Theorem \ref{near-parallel-k}, since it is a natural generalization of Theorem \ref{near-parallel}, with an only slightly more involved proof. 
As of now, we don't know how to use it to give  a better analysis of branching in \eref{ss-r}. 
It is, however, intuitive, and nice, that the approximation guarantee in \eref{sinak} gets better and better, as $k$ grows. 
}
Let us define 
\beq
\ba{rcl}
f(a,k)  & = &  2^{(k(n-k)+1)/4}/\norm{a}^{k/n} \\
g(a,k)  & = &  2^{k(n-1-k)/4}/\norm{a}^{(k-1)/n}.  
\ea
\eeq

\nin The successive version of Theorem \ref{near-parallel} is given below:
\bth \label{near-parallel-k}
\label{pa-1-k}
Let $a \in \zad{n} \,$ be a row-vector, with $\norm{a} \geq 2^{(n/2 +1)n}, \,$  $U$ a unimodular matrix such that 
the columns of 
$$
\bpx
a \\
I
\epx U 
$$
are LLL-reduced, and 
$P_k \,$ the (integral) submatrix of $U^{-1}$ consisting of 
the last $k$ rows. Furthermore, let 
$a(k)$ be the projection of $a$ onto the subspace spanned
by the rows of $P_k, \, r = a - a(k), \,$ and 
$$
\lambda_k := \norm{a(k)}/\mydet(P_k P_k^T)^{1/2}.
$$
Then 
\benum
\item \label{blah} $(\mydet(P_k P_k^T))^{1/2} (1+\norm{r}^{2})^{1/2} \leq \norm{a} f(a,k)$;
\item $\lambda_k \geq 1/f(a,k)$;
\item \label{sinak} $|\sin(a,a(k))| \leq \norm{r}/\lambda_k \leq 2 f(a,k)$.
\eenum
\enth
\pf{sketch} We will use the notation of Theorem \ref{near-parallel}. 
In its proof we simply change \eref{lnm1-1} (we copy the first expression for $\mydet L_n \,$ for easy reference) to 
\beq \label{lnm1-k}
\ba{rcl} 
\mydet L_{n} & = & \mydet L_{n}^\perp \, = \, (\norm{a}^2 + 1)^{1/2},  \\ 
\mydet L_{n-k} & = & \mydet L_{n-k}^\perp \, = \, (\mydet(P_k P_k^T))^{1/2} (1+\norm{r}^{2})^{1/2},
\end{array}
\eeq
and \eref{mydet} to 
\beq \label{mydet-k}
\ba{rcl}
\mydet \, L_{n-k} & \leq & 2^{k(n-k)/4} (\mydet L_n)^{1 - k/n}.
\end{array}
\eeq
Then substituting into \eref{mydet-k} from \eref{lnm1-k} gives 
\beq
\ba{rcl}
(\mydet(P_k P_k^T))^{1/2} (1+\norm{r}^{2})^{1/2} & \leq & 2^{(k(n-k))/4} (\sqrt{\norm{a}^2 + 1})^{1 - k/n}  \\
                             & \leq & 2^{(k(n-k)+1)/4}/\norm{a}^{k/n} \\ 
                             & = & \norm{a} f(a,k),
\end{array}
\eeq
with the second  inequality coming the lower bound on $\norm{a}$. This shows \eref{blah}, and the rest of the proof follows verbatim the proof 
of Theorem \ref{near-parallel}.
\qed

\nin Theorem \ref{near-parallel-null} also has a successive variant, which is 
\bth \label{near-parallel-null-k}
Suppose $\norm{a} \, \geq \, 2^{(n/2 +1)n}. \,$ Let 
$V$ be a matrix whose columns are an LLL-reduced basis of $\nlatt{a}$, $b$ an integral column vector with $ab=1$, $k \leq n-1 \,$ an integer,
and $P_k$ the (integral) submatrix of $(V,b)^{-1}$ consisting of the next-to-last $k \,$ rows. 

\nin Furthermore, let $a(k)$ be the projection of $a$ onto the subspace spanned
by the rows of $P_k, \, r = a - a(k), \,$ and 
$$
\lambda_k := \norm{a(k)}/\mydet(P_k P_k^T)^{1/2}.
$$
Then $r \neq 0, \,$ and 
\benum
\item \label{blah-null} \label{pa-null-1-1-k} $(\mydet(P_k P_k^T))^{1/2} \norm{r} \leq \norm{a} g(a,k)$;
\item \label{pa-null-1-2-k} $|\sin(a,a(k))| \leq \norm{r}/\lambda \leq 2 g(a,k)$.
\eenum
\enth
\pf{sketch} We will use the notation of Theorem \ref{near-parallel-null}. We need to replace \eref{lnm1-2} with 
\beq \label{lnm1-2-k}
\ba{rcl} 
\mydet L_{n-1} & = & \mydet L_{n-1}^\perp \, = \, \norm{a},  \\ 
\mydet L_{n-1-k} & = & \mydet L_{n-1-k}^\perp \, = \, (\mydet(P_k P_k^T))^{1/2} \norm{r}. 
\end{array}
\eeq
Theorem \ref{sublat-thm} with $n-1$ in place of $n, \,$ and $n-1-k$ in place of $\ell$ implies 
\beq \label{mydet-2-k}
\ba{rcl}
\mydet \, L_{n-1-k} & \leq & 2^{k(n-1-k)/4} (\mydet L_{n-1})^{1 - k/(n-1)}.
\end{array}
\eeq
Plugging the expressions for  $\mydet L_{n-1} \,$ and $\mydet L_{n-1-k} \,$ from \eref{lnm1-2-k} into \eref{mydet-2-k} gives
\beq
\ba{rcl}
(\mydet(P_k P_k^T))^{1/2} \norm{r} & \leq & 2^{k(n-1-k)/4} \norm{a}^{1 - k/(n-1)} \\
                                   & = & g(a,k) \norm{a},
\ea
\eeq
proving \eref{blah-null}. The rest of the proof is an almost verbatim copy of the corresponding proof in Theorem \ref{near-parallel-null}.
\qed

\co{

\subsection{Width vs integer width}
\label{sub-width}

The main goal of this paper was to prove a strong upper bound ???
}
\co{
\section{Conclusion} 

We proved upper bounds 
}

\nin{\bf Acknowledgement} We thank Don Coppersmith for his generous, and kind help on the $n=2$ case.
Thanks are due to Ravi Kannan for helpful discussions; to Laci Lov\'asz and Fritz Eisenbrand for discussions on the connection 
with diophantine approximation; and to  Jeff Lagarias and Andrew Odlyzko for pointing out reference \cite{CJLOSS92}. 

\bibliography{IP_Refs}

\begin{thebibliography}{10}

\bibitem{ABHLS00}
Karen Aardal, Robert~E. Bixby, Cor A.~J. Hurkens, Arjen~K. Lenstra, and Job~W.
  Smeltink.
\newblock Market split and basis reduction: Towards a solution of the
  {C}ornu{\'{e}}jols-{D}awande instances.
\newblock {\em INFORMS Journal on Computing}, 12(3):192--202, 2000.

\bibitem{AHL00}
Karen Aardal, Cor A.~J. Hurkens, and Arjen~K. Lenstra.
\newblock Solving a system of linear {D}iophantine equations with lower and
  upper bounds on the variables.
\newblock {\em Mathematics of Operations Research}, 25(3):427--442, 2000.

\bibitem{AL04}
Karen Aardal and Arjen~K. Lenstra.
\newblock Hard equality constrained integer knapsacks.
\newblock {\em Mathematics of Operations Research}, 29(3):724--738, 2004.

\bibitem{CRSS93}
William Cook, Thomas Rutherford, Herbert~E. Scarf, and David~F. Shallcross.
\newblock An implementation of the generalized basis reduction algorithm for
  integer programming.
\newblock {\em ORSA Journal on Computing}, 5(2):206--212, 1993.

\bibitem{CD98}
G{\'{e}}rard Cornu{\'{e}}jols and Milind Dawande.
\newblock A class of hard small 0--1 programs.
\newblock In {\em 6th Conference on Integer Programming and Combinatorial
  Optimization}, volume 1412 of {\em Lecture notes in Computer Science}, pages
  284--293. Springer-Verlag, 1998.

\bibitem{CJLOSS92}
M.~J. Coster, A.~Joux, B.~A. LaMacchia, A.~M. Odlyzko, C.~P. Schnorr, and
  J.~Stern.
\newblock Improved low-density subset sum algorithms.
\newblock {\em Computational Complexity}, 2:111--128, 1992.

\bibitem{EL05}
Friedrich Eisenbrand and S{\"{o}}ren Laue.
\newblock A linear algorithm for integer programming in the plane.
\newblock {\em Mathematical Programming}, 102(2):249--259, 2005.

\bibitem{FT87}
Andr\'as Frank and \'{E}va Tardos.
\newblock An application of simultaneous diophantine approximation in
  combinatorial optimization.
\newblock {\em Combinatorica}, 7(1):49--65, 1987.

\bibitem{FK89}
Merrick Furst and Ravi Kannan.
\newblock Succinct certificates for almost all subset sum problems.
\newblock {\em SIAM Journal on Computing}, 18:550 -- 558, 1989.

\bibitem{GZ02}
Liyan Gao and Yin Zhang.
\newblock Computational experience with lenstra's algorithm.
\newblock {\em Technical Report, Department of Computational and Applied
  Mathematics, Rice University}, 2002.

\bibitem{HN04}
Walfred Huyer and Arnold Neumaier.
\newblock Integral approximation of rays and verification of feasibility.
\newblock {\em Reliable Computing}, 10:195--207, 2004.

\bibitem{K87}
Ravi Kannan.
\newblock Minkowski's convex body theorem and integer programming.
\newblock {\em Mathematics of Operations Research}, 12(3):415--440, 1987.

\bibitem{KZ1873}
A.~Korkine and G.~Zolotarev.
\newblock Sur les formes quadratiques.
\newblock {\em Mathematische Annalen}, 6:366--389, 1873.

\bibitem{KP06}
Bala Krishnamoorthy and G{\'{a}}bor Pataki.
\newblock Column basis reduction and decomposable knapsack problems.
\newblock {\em Research Report 2006-07, Dept of Statistics and Operations
  Research, UNC-Chapel Hill, under review,
  http://www.optimization-online.org/DB\_{}HTML/2007/06/1701.html,
  http://arxiv.org/abs/0807.1317}, 2006.

\bibitem{L85}
Jeffrey~C. Lagarias.
\newblock The computational complexity of simultaneous diophantine
  approximation.
\newblock {\em SIAM J. Comput.}, 14:196--209, 1985.

\bibitem{LO85}
Jeffrey~C. Lagarias and Andrew~M. Odlyzko.
\newblock Solving low-density subset sum problems.
\newblock {\em Journal of ACM}, 32:229--246, 1985.

\bibitem{LLL82}
Arjen~K. Lenstra, Hendrik~W. {Lenstra, Jr.}, and L{\'{a}}szl{\'{o}}
  Lov{\'{a}}sz.
\newblock Factoring polynomials with rational coefficients.
\newblock {\em Mathematische Annalen}, 261:515--534, 1982.

\bibitem{L83}
Hendrik~W. {Lenstra, Jr.}
\newblock Integer programming with a fixed number of variables.
\newblock {\em Mathematics of Operations Research}, 8:538--548, 1983.

\bibitem{LS92}
L{\'{a}}szl{\'{o}} Lov{\'{a}}sz and Herbert~E. Scarf.
\newblock The generalized basis reduction algorithm.
\newblock {\em Mathematics of Operations Research}, 17:751--764, 1992.

\bibitem{M03}
Jacques Martinet.
\newblock {\em Perfect Lattices in Euclidean Spaces}.
\newblock Springer-Verlag, Berlin, 2003.

\bibitem{ML04}
Sanjay Mehrotra and Zhifeng Li.
\newblock On generalized branching methods for mixed integer programming.
\newblock {\em Research Report, Department of Industrial Engineering,
  Northwestern University}, 2004.

\bibitem{PT08}
G\'abor Pataki and Mustafa Tural.
\newblock On sublattice determinants in reduced bases.
\newblock {\em Technical Report 2008-02, Dept of Statistics and Operations
  Research, UNC Chapel Hill, under review,
  http://www.optimization-online.org/DB\_HTML/2008/04/1960.html,
  http://arxiv.org/abs/0804.4014}.

\bibitem{S86}
Alexander Schrijver.
\newblock {\em Theory of Linear and Integer Programming}.
\newblock Wiley, Chichester, United Kingdom, 1986.

\end{thebibliography}

\end{document}